\documentclass[12pt,draftcls,onecolumn]{IEEEtran}

\usepackage{amsfonts}
\usepackage{amssymb}
\usepackage{amsmath}
\usepackage{mathptmx}
\usepackage{cancel}
\usepackage{pstricks}
\usepackage{psfrag}
\usepackage{mathrsfs}
\usepackage{graphicx}
\usepackage{pgf,tikz}
\usepackage{algorithm}
\usepackage{algpseudocode}
\usepackage{soul}
\usetikzlibrary{arrows}
\def\be{\begin{equation}}
\def\ee{\end{equation}}
\def\bea{\begin{eqnarray}}
\def\eea{\end{eqnarray}}
\def\beann{\begin{eqnarray*}}
\def\eeann{\end{eqnarray*}}
\def\spacingset#1{\def\baselinestretch{#1}\small\normalsize}
\spacingset{1.35}

\newtheorem{theorem}{Theorem}

\newtheorem{definition}{Definition}

\def\be{\begin{equation}}
\def\ee{\end{equation}}
\def\bea{\begin{eqnarray}}
\def\eea{\end{eqnarray}}
\def\beann{\begin{eqnarray*}}
\def\eeann{\end{eqnarray*}}

\def\proof{\noindent{\bf{\em Proof:}\ \ }}
\def\QED{\mbox{\rule[0pt]{1.5ex}{1.5ex}}}
\def\endproof{\hspace*{\fill}~\QED\par\endtrivlist\unskip}

\newcommand{\beq}{\begin{equation}}
\newcommand{\eeq}{\end{equation}}

\def\bmat{\left[ \begin{array}}
\def\emat{\end{array} \right]}

\def\bmat{\left[ \begin{array}}
\def\emat{\end{array} \right]}
\def\bsmat{\left[ \begin{smallmatrix}}
\def\esmat{\end{smallmatrix} \right]}

\def\i{\mathfrak{i}}

\newcommand{\ba}{\begin{array}}
\newcommand{\ea}{\end{array}}

\begin{document}
\title{\LARGE{A study on Vandermonde-like polynomial matrices}}

\author{Augusto Ferrante,  Fabrizio Padula, and Lorenzo Ntogramatzidis

 \thanks{Augusto Ferrante is with the Dipartimento di Ingegneria dell'Informazione, Universit\`a di Padova, 
	via Gradenigo, 6/B -- I-35131 Padova, Italy. E-mail: {\tt augusto@dei.unipd.it. }}
	 \thanks{Fabrizio Padula and L. Ntogramatzidis are with the Department of Mathematics and
Statistics, Curtin University, Perth,
Australia. E-mail: {\tt \{Fabrizio.Padula,L.Ntogramatzidis\}@curtin.edu.au. }}



}


\maketitle

\vspace{-1cm}

\IEEEpeerreviewmaketitle

\begin{abstract}
A new class of structured matrices is presented and a closed form formula for their determinant is established. This formula has strong connections with the one for  Vandermonde matrices. 
\end{abstract}


\section{Introduction}
\label{secintro}

In this brief note we introduce a new class of structured matrices that we call {\em controllability intermixing matrices}, or {\em CI-matrices},
because of their role in a problem of dynamics assignment in linear system theory, see 
\cite{PFN-EigenstructureAssignmentinGeometricControl-arxiv}, which is our motivation for studying these matrices.
For this class of matrices, we establish a formula much in the same vein of the one holding  for the celebrated Vandermonde matrices, \cite[page 37]{Horn-Johnson}.
This formula allows to compute the determinant  of  CI-matrices  in closed-form and, as shown in 
\cite{PFN-EigenstructureAssignmentinGeometricControl-arxiv}, turns out to be the key result to prove
the invariance in the dimension of the projection of the null-spaces of the controllability and of the Rosenbrock matrix pencils onto the state space of a linear time-invariant system.



\subsection{Controllability intermixing matrices}
Let $n\in \mathbb N$ and $\mathscr{M}:=\{\mu_1,\mu_2,\dots,\mu_n\}$ be a set of indeterminates.
Let $\mathscr{M}_k:=\{\mu_i\in\mathscr{M}:\ i\neq k  \}$ and 
let $\mathscr{I}_k:=\{1,2,\dots,k-1,k+1,\dots,n\}$ be the set of indices of $\mathscr{M}_k$.
\begin{definition}
We say that $M$ is a {\em controllability intermixing matrix}, or {\em CI-matrix} of size $n$ if it is an
$n\times n$ matrix, polynomial in $\mu_1,\mu_2,\dots,\mu_n$, of dimensions $n\times n$, whose entries  are
\beq
M_{h,k}=\left\{\begin{array}{cl}
{\displaystyle \sum_{{\scriptsize\begin{array}{c}i_1,i_2,\dots, i_{n-h} \in\mathscr{I}_k,\\[-1mm] i_1<i_2<\dots<i_{n-h}\end{array}}}} \mu_{i_1}\mu_{i_2}\dots \mu_{i_{n-h}}& {\rm \ \ \  if\ } h<n\\
 1 & {\rm\ \ \  if\ } h=n.\\
 \end{array}\right.
\eeq
\end{definition}

A {\em CI-matrix} of size $n$ has thus the following structure 
\beq\label{defmatricetipovdm}
M=\bmat{llcl}
\ \ \ \ \ {\displaystyle \prod_{i\in\mathscr{I}_1}\mu_{i} } & \ \ \ \ \ {\displaystyle \prod_{i\in\mathscr{I}_2}\mu_{i} }&\dots
&\ \ \ \ \ {\displaystyle \prod_{i\in\mathscr{I}_n}\mu_{i} } \\
\ \ \ \ \ \ \   \vdots&\ \ \ \ \ \ \  \vdots&\vdots&\ \ \ \ \ \ \  \vdots\\
{\displaystyle \sum_{{\scriptsize\begin{array}{c}i_1,i_2,i_3 \in\mathscr{I}_1,\\[-1mm]  i_1<i_2<i_3\end{array}}}} \mu_{i_1}\mu_{i_2}\mu_{i_3} & {\displaystyle \sum_{{\scriptsize\begin{array}{c}i_1,i_2,i_3 \in\mathscr{I}_2,\\[-1mm]  i_1<i_2<i_3\end{array}}}} \mu_{i_1}\mu_{i_2}\mu_{i_3} & \dots &{\displaystyle \sum_{{\scriptsize\begin{array}{c}i_1,i_2,i_3 \in\mathscr{I}_n,\\[-1mm]  i_1<i_2<i_3\end{array}}}} \mu_{i_1}\mu_{i_2}\mu_{i_3}\\
\ {\displaystyle \sum_{{\scriptsize\begin{array}{c} i_1,i_2\in\mathscr{I}_1,\\[-1mm]  \i_1<i_2\end{array}}}} \mu_{i_1}\mu_{i_2} & \ {\displaystyle \sum_{{\scriptsize\begin{array}{c} i_1,i_2\in\mathscr{I}_2,\\[-1mm]  \i_1<i_2\end{array}}}} \mu_{i_1}\mu_{i_2}& \dots & \ {\displaystyle \sum_{{\scriptsize\begin{array}{c} i_1,i_2\in\mathscr{I}_n,\\[-1mm]  \i_1<i_2\end{array}}}} \mu_{i_1}\mu_{i_2}\\
\ \ \ \ \ {\displaystyle\sum_{i\in\mathscr{I}_1}} \mu_i &\ \ \ \ \ {\displaystyle\sum_{i\in\mathscr{I}_2}} \mu_i &\dots &
\ \ \ \ \ {\displaystyle\sum_{i\in\mathscr{I}_n}} \mu_i\\
\ \ \ \ \ \ \ 1 &\ \ \ \ \ \ \ 1 &\dots &\ \ \ \ \ \ \ 1
\emat.
\eeq
To better understand the structure of $M$, we observe that the entries (from bottom to top) of its $k$-th column are: $1$, the sum of all the indeterminates of $\mathscr{M}_k$, the sum of the products of all the possible unordered pairs of indeterminates of $\mathscr{M}_k$, the sum of the products of all the possible unordered triples of indeterminates of $\mathscr{M}_k$ and so on up to 
the sum of the products of all the possible unordered $(n-1)$-tuples of indeterminates of $\mathscr{M}_k$ (notice that there is just one $(n-1)$-tuple of indeterminates of $\mathscr{M}_k$, i.e. the $(n-1)$-tuple containing all the elements of $\mathscr{M}_k$, so that the sum in the first entry of each column is over just one element and for this reason this sum is missing in formula \eqref{defmatricetipovdm}).
Again to help intuition, we write below the full matrix for $n=4$:
{\small
\beann
\bmat{ccccc}
 \mu_{\scriptscriptstyle 2} \mu_{\scriptscriptstyle 3} \mu_{\scriptscriptstyle 4}  & \mu_{\scriptscriptstyle 1} \mu_{\scriptscriptstyle 3} \mu_{\scriptscriptstyle 4}   &  \mu_{\scriptscriptstyle 1} \mu_{\scriptscriptstyle 2} \mu_{\scriptscriptstyle 4}   &  \mu_{\scriptscriptstyle 1} \mu_{\scriptscriptstyle 2} \mu_{\scriptscriptstyle 3}   \\
 (\mu_{\scriptscriptstyle 2} \mu_{\scriptscriptstyle 3} + \mu_{\scriptscriptstyle 3} \mu_{\scriptscriptstyle 4}+ \mu_{\scriptscriptstyle 2} \mu_{\scriptscriptstyle 4})  &   (\mu_{\scriptscriptstyle 1} \mu_{\scriptscriptstyle 3} + \mu_{\scriptscriptstyle 3} \mu_{\scriptscriptstyle 4}+ \mu_{\scriptscriptstyle 1} \mu_{\scriptscriptstyle 4})  &  (\mu_{\scriptscriptstyle 1} \mu_{\scriptscriptstyle 2} + \mu_{\scriptscriptstyle 2} \mu_{\scriptscriptstyle 4}+ \mu_{\scriptscriptstyle 1} \mu_{\scriptscriptstyle 4})  &   (\mu_{\scriptscriptstyle 1} \mu_{\scriptscriptstyle 2} + \mu_{\scriptscriptstyle 2} \mu_{\scriptscriptstyle 3}+ \mu_{\scriptscriptstyle 1} \mu_{\scriptscriptstyle 3})  \\
 (\mu_{\scriptscriptstyle 2}+\mu_{\scriptscriptstyle 3}+\mu_{\scriptscriptstyle 4})  &   (\mu_{\scriptscriptstyle 1}+\mu_{\scriptscriptstyle 3}+\mu_{\scriptscriptstyle 4})  &   (\mu_{\scriptscriptstyle 1}+\mu_{\scriptscriptstyle 2}+\mu_{\scriptscriptstyle 4})  &   (\mu_{\scriptscriptstyle 1}+\mu_{\scriptscriptstyle 2}+\mu_{\scriptscriptstyle 3})  \\
1 &1 & 1 & 1 \emat.
\eeann 
}

\section{Determinant of CI-matrices}
We now establish a formula for the computation of the determinant of CI-matrices.
\begin{theorem}
The determinant of the CI-matrix $M$ defined by (\ref{defmatricetipovdm}) is  given by
\beq
\det(M)=\prod_{1\leq i<j\leq n} (\mu_{j}-\mu_{i}).
\eeq
\end{theorem}
\proof
All the entries of $M$ are polynomials in the indeterminates $\mu_i$.
Moreover, all the entries in the $j$-th row of $M$ have the same degree: this degree is 
indeed $n-j$.
Hence, in view of Leibniz formula, all the monomials of the polynomial $\det(M)$ have   degree equal to 
\beq
\sum_{j=1}^n (n-j) = n(n-1)/2;
\eeq
in other words, the determinant $\det(M)$  is a homogeneous polynomial of  degree $n(n-1)/2$.
Now observe that if $\mu_{i}=\mu_{j}$ the columns $m_i$ and $m_j$ of $M$ are equal by construction, so that $\det(M)=0$. As a consequence, for all $1\leq i<j\leq n$, the polynomials $\mu_{j}-\mu_{i}$ 
divide $\det(M)$. Moreover, these polynomials are coprime, so that also  their product
divides $\det(M)$; in other words
\beq
\det(M)=N\prod_{1\leq i<j\leq n} (\mu_{j}-\mu_{i}),
\eeq
with $N$ being also a polynomial. But 
\beq
\deg\left[\prod_{1\leq i<j\leq n} (\mu_{j}-\mu_{i})\right]=n(n-1)/2=\deg[\det(M)],
\eeq
so that clearly $N$ is a constant.
It remain to show that $N=1$.
This will be show inductively on $n$.
To this end we denote by $M_n$ the matrix of this class associated with $n$ indeterminates
and by $N_n$ the corresponding constant $N$. 
The base of induction is obvious as $M_1=1$, so that $N_1=1$.
Assume now that $N_n=1$ and consider
the matrix $M_{n+1}$ associated with the $n+1$ indeterminates $\mu_1,\mu_2,\dots,\mu_{n+1}$.
We know that 
\beq
\det(M_{n+1})=N_{n+1} \prod_{1\leq i<j\leq n+1} (\mu_{j}-\mu_{i})
\eeq
and we want to show that $N_{n+1}=1$.
For $\mu_{1}=0$, the latter formula reads as
\beq\label{primaformadetmhp1}
\det(M_{n+1})=N_{n+1} \prod_{2\leq i \leq n+1} \mu_{i}
\prod_{2\leq i<j\leq n+1} (\mu_{j}-\mu_{i}).
\eeq
On the other hand, for $\mu_{1}=0$, we see by inspection  that
\beq
M_{n+1}=\bmat{cc}{\displaystyle \prod_{1\leq i<j\leq n+1}} \mu_{i}& 0\\
\ast& M_{n}\emat,
\eeq
where the matrix $M_{n}$ is a CI-matrix of size $n$ associated with the $n$ indeterminates $\mu_2,\dots,\mu_{n+1}$.
Therefore, we have
\beq
\det(M_{n+1})=\left[\prod_{2\leq i \leq n+1} \mu_{i}\right]
\det(M_{n}),
\eeq
and by the inductive assumption we have
\beq
\det(M_{n+1})=\prod_{2\leq i \leq n+1} \mu_{i}
\prod_{2\leq i<j\leq n+1} (\mu_{j}-\mu_{i}).
\eeq
Comparing the latter with (\ref{primaformadetmhp1}), we get $N_{n+1}=1$ as desired.
\endproof



\begin{thebibliography}{1}
	



\bibitem{Horn-Johnson}
R. A. Horn, and C. R. Johnson.
\newblock {\em Matrix Analysis}.
\newblock Cambridge University Press, 2013.

	

\bibitem{PFN-EigenstructureAssignmentinGeometricControl-arxiv}
F. Padula, A. Ferrante, and L.~Ntogramatzidis.
\newblock 	New geometric results in eigenstructure assignment, {\tt https://arxiv.org/abs/1910.10867}, 2019.
	
	\end{thebibliography}
\end{document}